\documentclass[10pt,letterpaper]{article}
\usepackage[top=0.85in,left=1.25in,footskip=0.75in]{geometry}
\usepackage{graphicx}
\graphicspath{{graphics/}{moregraphics/}}
\usepackage{xcolor}
\usepackage{epstopdf}
\usepackage[T1]{fontenc}
\usepackage[utf8]{inputenc}

\everymath{\displaystyle}
\usepackage{cite}
\usepackage[cmex10]{amsmath}
\usepackage{amssymb,amsfonts}
\usepackage{array}
\usepackage{subfigure}
\usepackage[square,sort,comma,numbers]{natbib}
\usepackage{algorithmic}

\DeclareGraphicsExtensions{.pdf,.PDF,.jpeg,.JPEG,.jpg,.JPEG,.png,.PNG}

\usepackage{booktabs}
\usepackage{textcomp}
\usepackage[spaces,hyphens]{url}

\def\BibTeX{{\rm B\kern-.05em{\sc i\kern-.025em b}\kern-.08em
    T\kern-.1667em\lower.7ex\hbox{E}\kern-.125emX}}

\begin{document}
\begin{flushleft}
{\Large
\textbf\newline{A Study of a Loss System with Priorities} 
}
\newline
\\
Hang Yang\textsuperscript{1},
Jing Fu\textsuperscript{2},
Jingjin Wu\textsuperscript{3,4,*},
Moshe Zukerman\textsuperscript{1}

\bigskip
\textbf{1} Department of Electrical Engineering, City University of Hong Kong, Kowloon, Hong Kong SAR, P. R. China
\\
\textbf{2} School of Engineering, RMIT University, Melbourne, Victoria 3000, Australia
\\
\textbf{3} Department of Statistics and Data Science, BNU-HKBU United International College, Zhuhai, Guangdong, 519087, P. R. China
\\
\textbf{4} Guangdong Provincial Key Laboratory of Interdisciplinary Research and Application for Data Science, Guangdong, 519087, P. R. China
\bigskip

%
%

* Corresponding author, E-mail: 
jj.wu@ieee.org 

\end{flushleft}




\section*{Abstract}

 The Erlang loss formula, also known as the Erlang~B formula, has been known for over a century and has been used in a wide range of applications, from telephony to hospital intensive care unit management. It provides the blocking probability of arriving customers to a loss system involving a finite number of servers without a waiting room. 
  Because of the need to introduce
priorities in many services, an extension of the Erlang B formula to the case
of a loss system with preemptive priority is valuable and essential.
  This paper analytically establishes the consistency between the global balance (steady state) equations for a loss system with preemptive priorities and a known result obtained using traffic loss arguments for the same problem. This paper, for the first time, derives this known result directly from the global balance equations based on the relevant multidimensional Markov chain. The paper also addresses the question of whether or not the well-known insensitivity property of the Erlang loss system is also applicable to the case of a loss system with preemptive priorities, provides explanations, and demonstrates through simulations that, except for the blocking probability of the highest priority customers, the blocking probabilities of the other customers are sensitive to the service time distributions and that a larger service time variance leads to a lower blocking probability of the lower priority traffic.

\textbf{Keywords:} Loss system, Erlang B system, preemptive priorities, multidimensional Markov chain, insensitivity


\section{Introduction}\label{sec:intro1}
In many applications, there are service resources abstracted as {\it servers}, and incoming customers that are either served by available service resources, blocked and cleared from the system, or overflowed to another service system.  Such systems are called {\it loss systems} and are characterized by a potential lack of waiting room. Examples of service resources in such a system include telephone circuits, wireless channels, optical wavelengths, and hospital beds of intensive care units. In this paper, we use the term ``customers'' to refer to a range of requests seeking appropriate services, such as phone calls, requests for wireless or wavelength channels, as well as actual customers, patients, clients, or travelers.   

A loss system where arriving customers follow a Poisson process and their service times are independent and exponentially distributed is denoted by M/M/$k$/$k$ (Kendall notation~\citep{Kendall53}) where the first $k$ represents the number of  servers, and the second $k$ is the maximal number of service requests allowed in the system or the number of buffer places including the buffer places at the servers. Throughout the paper, the terms ``service time'' and ``holding time'' are considered synonymous and will be used interchangeably.  

Applications of loss systems in general, and M/M/$k$/$k$ in particular, including critical systems such as Intensive Care Units (ICU), where the availability of a service resource may have life-and-death consequences~\citep{Litvak08}. In other applications, such as those in the areas of computers and communications systems and networks, excessive blocking events have adverse consequences on customers' Quality of Service (QoS).  Erlang Loss Formula (also known as Erlang B Formula)  is a century-old solution for the blocking probability of the arriving customers which has been widely applied in many areas, such as 
hospital resource allocation, \citep{deBruin10,Bekker10,Bekker17,Andersen17,Andersen19,Vanberkel11}, 
telephony \citep{Adams24,bray1986telecommunications,Kelly86}, 
mobile networks \citep{Eklundh86,Hong86,
Everitt89,zhai2016spectrum,hiew2000,chen2014extra,WU2017},
video on demand~\citep{De_Giovanni94,Li96}, 
call centers \citep{Vakilinia2015},  parallel computing systems~\citep{Hyytia2019}, self-driving cars \citep{hampshire2020beyond} and
optical networks \citep{Washington2004,Li19,Wang14,ZALESKY2007}.
Over the years, there have been publications that considered extensions of loss systems beyond Erlang B; they include \citep{brandwajn2017multi,Glabowski21,klimenoklack,klimenok2020priority,Moscholios19,Li2019}. 



The Erlang Loss Formula is also known to be insensitive to the service time distribution beyond its mean  \citep{Sevastyanov57}. In other words, the formula gives exact blocking probability results regardless of whether or not the service time is exponentially distributed, so as long as the mean service time is known, its actual distribution does not affect the resulted blocking probabilities, and the formula provides exact results for any service time distribution.  This insensitivity property is key to the broad applicability of the Erlang Loss Formula in many practical cases where the service time does not follow an exponential distribution. 

In various applications, services are differentiated according to various criteria. For example, in a hospital triage process, patients are differentiated according to the severity of their conditions \citep{Ding19}. In other service systems, such as packet switching networks~\citep{Li2023} and downlink transmission in 5G mobile networks for ultra-reliable low latency communication (URLLC) traffics~\citep{Kim2020}, customer requests may be differentiated according to the cost and/or their QoS requirements. In such cases, implementation of preemptive priorities is an option~\citep{Barros08,Pal19,Palit20,Saad18}. 


Another application of preemptive priorities is to approximate blocking probabilities in overflow loss systems by constructing a surrogate system with preemptive priorities~\citep{WU2017, Wong2007, Wu2020, CHAN2016}. However, in this paper, we focus on deriving exact evaluations of blocking probabilities from the global balance equations  based on multidimensional Markov Chains and thus have a fundamentally different objective compared to the studies similar to~\citep{WU2017, Wong2007, Wu2020,CHAN2016}.

We consider a Markovian loss system with a finite number of servers where arriving customers are classified into a finite number of preemptive-priority classes. The arrival process of each class of customers follows a Poisson process, and all the service times are independent and exponentially distributed with a given parameter. In this Markovian loss system, if a customer arrives and all the servers are busy, the arriving customer may preempt a customer in the system that is being served if the arriving customer has higher priority than the preempted customer. The problem we consider is to find the blocking probability of customers from each of the preemptive priority classes.  

The first to consider this problem was Katzschner \citep{Katzschner70} in 1970. He considered the global balance equations of this loss system with preemptive priorities. He then derived the Laplace transform of the joint steady-state distribution of the number of customers of each class for the case of two classes from which the blocking probabilities can be obtained. Finally, he explained how this solution could be extended to the general case of a finite number of customer classes. Vu and Zukerman \citep{V02}, who were not aware of \citep{Katzschner70}, provided a simple procedure based on rate arguments to derive the blocking probabilities for all priority classes in the context of optical burst switching (OBS) application. Then, Yang and Stol \citep{Yang14} (who were not aware of  \citep{Katzschner70} and \citep{V02}) obtained the blocking probabilities numerically from the global balance steady-state equations. None of the previous papers that dealt with this loss system with preemptive priorities established the equivalence between the global balance steady-state equations and the simple solution of \citep{V02} analytically. In this paper, we do precisely that; we start from the global balance equations, and we rigorously prove the solution of \citep{V02}. 

Another important contribution of this paper is to address the fundamental question of whether or not the system we consider, i.e., a loss system with preemptive priorities, is insensitive to service time distribution beyond its mean (i.e. to the shape of the service time distribution) of the customers of various service classes. As mentioned above, it is well known that the Erlang B formula is insensitive to the shape of the service time distribution \citep{Sevastyanov57}. The question of whether or not this insensitivity property is also applicable to the blocking probabilities of the various classes of customers in a loss system with multiple priorities has not been properly addressed for over half a century since this problem was first introduced by Katzschner \citep{Katzschner70} in 1970.  
In this paper, we explain and demonstrate by simulations that for the highest preemptive priority that is not affected by lower priority traffic the insensitivity of the Erlang B system applies. However, for lower preemptive priorities, insensitivity does not hold, and that higher variance of the service time distribution leads to lower blocking probabilities of the lower priority traffic. This has not been demonstrated in the earlier publications on this problem, including \citep{Katzschner70,V02,Yang14}. In fact, in \citep{V02}, an incorrect comment was made that implies the insensitivity property for all priority classes. 

The remainder of the paper is organized as follows. In Section \ref{V2}, we describe the model and provide the solution of  \citep{V02} for the blocking probability of the customers of the different priority classes. In Section \ref{sec:multi}, we provide a multi-dimensional Markov-chain analysis for our problem and establish consistency with the results of \citep{V02} presented in Section \ref{V2}. In Section~\ref{sec:sens}, we explain and demonstrate by simulations that except for the blocking probability of the highest priority customers,  the blocking probabilities of the other customers are sensitive to the holding time distributions. 
Finally, the conclusions are drawn in Section~\ref{sec:conc}.


\section{Model Description and the Solution of Vu and Zukerman}\label{V2}

This paper considers the fundamental problem of obtaining the blocking probability in a loss system where the customers are classified into $p$ preemptive priority classes and the arrival process of priority $i$ customers  (for $i=1,2,3, \ldots, p$) follows a Poisson process with parameter $\lambda_i$. The service time of the customers of all priorities is exponentially distributed with parameter $\mu$. Accordingly, the offered traffic of the customers of priority $i$  is \begin{equation}
\label{Aifixedmu}
A_i = \frac{\lambda_i}{\mu}, ~~ i=1,2,3, \ldots, p.\end{equation}
Priority 1 is the highest priority, and priority $p$ is the lowest. If $i>j$, then customers of priority $i$ have lower priority than customers of priority $j$. In this case, an arriving priority $j$ customer, at the time of its arrival, may preempt a priority $i$ customer already being served.

This section provides the solution of Vu and Zukerman \cite{V02} based on rate arguments for the blocking probabilities of customers of each priority class. This solution is also described in \cite{Zukerman20}.



Let $P_{\text{b}}(i)$ denote the blocking probability of the customer class with priority $i$. Because priority 1 customers may preempt lower priority customers, they can access the loss system as if these lower priority customers do not exist.  Therefore, for $i=1$, we obtain
\begin{equation}P_{\text{b}}(1)= E_k (A_1),\end{equation}
where $E_k(A)$ denotes the blocking probability obtained by the Erlang B formula for an M/M/$k$/$k$ loss system with $k$ servers and offered traffic $A$. Specifically,

\begin{equation}
E_k(A) = \dfrac{A^k/k!}{\displaystyle\sum\limits_{m=0}^k\dfrac{A^m}{m!}}.
\end{equation}

To obtain $P_{\text{b}}(i)$ for $i =  2, 3, \ldots, p$, the first step is to observe \citep{V02,Zukerman20} that the blocking probability of the total traffic from priority $i$ and higher, which is the traffic generated by customers of priorities $1, 2, \ldots, i$, is equal to:
$$E_k\left(\displaystyle\sum_{j=1}^i A_j \right).$$

To explain this observation, we need to consider the fact that the overall blocking probability in the case of a loss system with a given number of servers and with a strict priority regime among multiple classes of traffic streams all having the same mean service time is equal to the blocking probability of a traditional loss system (with a single traffic class) with the same number of servers and the same mean service time as in the case of the multiple traffic classes.  This fact holds only because, in both cases, the service time distribution is exponential, and because this distribution is memoryless, a higher priority call that arrives when all servers are busy, which will be blocked in a single class loss system, will preempt a low property priority call in a multiple class loss system and will have the same remaining service time as that of the preempted low priority call. As for both calls, the remaining service time has an exponential distribution with mean $1/\mu$. Accordingly, in both systems, one blocked call will be recorded, and the remaining service times are also equal.


Henceforth, we will call this justification the {\it the displacement/replacement argument}, recalling the term ``displacing priorities'' of \citep{Katzschner70}.  Notice that using the displacement/replacement argument ignores the potential discrepancy in the priority system by replacing a preempted low-priority call with a high-priority call. 

However, the above-described method does not apply to cases where the service times are not exponentially distributed. 
In other words, except for the highest priority traffic, blocking probabilities results are not insensitive to the shape of holding time distribution (i.e., to the distribution beyond its mean). 
Unlike the conventional M/M/$k$/$k$ system, the blocking probabilities of the preempted customers with lower priorities are sensitive to the shape of their holding time distributions.
This will be further discussed, explained, and demonstrated by simulations in Section \ref{sec:sens}. This will correct \citep{V02} which implies that the insensitivity property applies to all priority classes. 


Next, we observe that the priority $i$ lost traffic, $i=2,3, \ldots, p$, denoted $A_L(i)$ is given by the lost traffic of priorities up to $i$ minus the lost traffic of priorities up to $i-1$, namely,
\begin{equation}
\label{lossdiff}
\begin{split}A_L(i) & = \left(\displaystyle\sum_{j=1}^i A_j \right)E_k \left(\displaystyle\sum_{j=1}^i A_j \right)  \\ & \quad  - \left(\displaystyle\sum_{j=1}^{i-1} A_j \right)E_k \left(\displaystyle\sum_{j=1}^{i-1} A_j \right).
\end{split} \end{equation}

Therefore, the $P_{\text{b}}(i)$ values for $2 \le i \le p$, are obtainable as the ratio of the priority $i$ lost traffic to the priority $i$ offered traffic,  that is,
\begin{equation}
\label{bp}
P_{\text{b}}(i)=\frac{A_L(i)}{A_i}.
\end{equation}

\section{Multi-dimensional Markov-chain Analysis}\label{sec:multi}
In this section, we provide a multi-dimensional Markov-chain analysis of our problem of a loss system with priorities. We start with the simpler case of 
two preemptive priorities. Then, we extend the analysis to the 
general case of $k$ servers and $p$ preemptive priorities.  
\subsection{The case of two preemptive priorities}\label{V2p}
We consider here the case of two classes of customers with two corresponding preemptive priorities. A similar analysis was also provided in \citep{Yang14}, but we include it here for completeness. 
The class $i$ customers are assumed to arrive following a Poisson process with parameter $\lambda_i$, $i=1, 2$. We consider a $k$-server loss system where customers are served by $k$ servers, and there is no waiting room.
We assume that the service times of all the customers are exponentially distributed with the parameter $\mu$. The class 1 customers have preemptive priority over the class 2 customers. That is, an arriving class 1 customer that finds all $k$ servers busy may preempt a class 2 customer in service and is served instead of the preempted class 2 customer. Accordingly, a class 1 (higher priority) customer will only be blocked if, upon its arrival, all the $k$ servers are busy serving class 1 customers, while a class 2 customer will be blocked on arrival if all $k$ servers are busy serving any of the customer classes, and may also be preempted during its service. 
Therefore, the offered traffic by class 1 customers is given by
\begin{equation} \label{eq1}
  \begin{split}
    A_1=\lambda_1/\mu.
  \end{split}
\end{equation}
and the offered traffic by class 2 customers is
\begin{equation} \label{eq2}
  \begin{split}
    A_2=\lambda_2/\mu.
  \end{split}
\end{equation}
Then, the total offered traffic is
\begin{equation} \label{eq3_2}
  \begin{split}
    A=A_1+A_2.
  \end{split}
\end{equation}

Let $(i,j)$ be the system's state, where $i$ represents the number of busy servers for class 1 customers and $j$ is the number of busy servers for class 2 customers. Thus, in state $(i,j)$,~ $i=0, 1, 2, \ldots, k$,~ $j =0, 1, 2, \ldots, k$, we must have 
\begin{equation} \label{eq4_2}
    0 \leq  i+j \leq k. 
\end{equation}
 Let $\pi_{i,j}$ be the steady-state probability of the system being in state $(i,j)$. Accordingly, 
\begin{equation} \label{eq5}
  \begin{split}
   \displaystyle\sum_{i=0}^k\sum_{j=0}^{k-i}\pi_{i,j}=1.
  \end{split}
\end{equation}

Recall that $P_{\text{b}}(1)$ represents the blocking probability of priority 1 customers (class 1 customers), and $P_{\text{b}}(2)$ is the blocking probability of priority 2 customers (class 2 customers). Note that $P_{\text{b}}(2)$ is the probability that an arriving priority 2 customer is either blocked on arrival or preempted by a priority 1 customer after it was admitted to service.  
Since priority 1 customers will only be blocked when all the $k$ servers are busy, and there are no priority 2 customers in service, we have  
\begin{equation} \label{eq6}
  \begin{split}
P_{\text{b}}(1)=\pi_{k,0}.
  \end{split}
\end{equation}

Let $P_{\text{boa}}(2)$ denote the probability that an arriving priority 2 customer finds all the servers busy and is blocked on arrival (boa), and let $P_{\text{pre}}(2)$ be 
the probability that an arbitrary priority 2 customer is lost due to the preemption of its service by an arrival of a priority 1 customer.


To obtain $P_{\text{boa}}(2)$, we observe that since the arriving priority 2 customers follow a Poisson process, and since an arriving priority 2 customer will be blocked on arrival if and only if all $k$ servers are busy, we have
\begin{equation} \label{eq7}
  \begin{split}
P_{\text{boa}}(2)=\displaystyle\sum_{i=0}^{k}\pi_{i,k-i}.
  \end{split}
\end{equation} 

An alternative approach to assert (\ref{eq7}) is to consider an arbitrarily long period of time $L$. Then, the mean number of priority 2 arrivals during $L$ is $\lambda_2 L$, and the mean number of priority 2 customers that are blocked on arrival during period $L$ is given by $\lambda_2 L \displaystyle\sum_{i=0}^{k}\pi_{i,k-i}$. Accordingly, the ratio of the latter to the former gives the proportion of priority 2 customers that are blocked on arrivals. 

To obtain $P_{\text{pre}}(2)$, we again consider an arbitrary long period of time $L$, 
and we already know that the mean number of priority 2 customers arriving during period $L$, is equal to $\lambda_2 L$ out of which  $P_{\text{boa}}(2)\lambda_2 L$ were blocked on arrival. We also observe that the mean number of priority  2 customers that are preempted during $L$ is equal to the number of priority 1 arrivals during the periods of times in $L$ that all servers are busy and there are some priority 2 customers being served because each one of these priority 1 arrivals will preempt a priority 2 customer that is being served. This mean is given by $$ \lambda_1L\displaystyle\sum_{i=0}^{k-1}\pi_{i,k-i}.$$
Then, $P_{\text{pre}}(2)$ is given by the ratio of the latter to the 
mean number of priority 2 customers that arrive during $L$ but are not blocked on arrival, so we obtain 
\begin{equation} \label{eq8}
  \begin{split}
P_{\text{pre}}(2) = \frac{\lambda_1}{
 \lambda_2 (1-P_{\text{boa}} (2))}\times\displaystyle\sum_{i=0}^{k-1}\pi_{i,k-i}.
  \end{split}
\end{equation} 


Having obtained $P_{\text{pre}}(2)$ and $P_{\text{boa}}(2)$, the overall blocking probability of a priority 2 customer is the probability that it is blocked on arrival ($P_{\text{boa}}(2)$), plus the probability that it is not blocked on arrival but lost due to the preemption of its service by an arrival of a priority 1 customer ($(1-P_{\text{boa}}(2))P_{\text{pre}}(2)$). Thus, 

\begin{equation} \label{eq9}
  \begin{split}
P_{\text{b}}(2) = P_{\text{boa}}(2) + (1-P_{\text{boa}}(2))P_{\text{pre}}(2) 
  \end{split}
\end{equation}

\begin{figure}[htpb]
\centering
\includegraphics[width=\columnwidth,height=\textwidth,keepaspectratio]{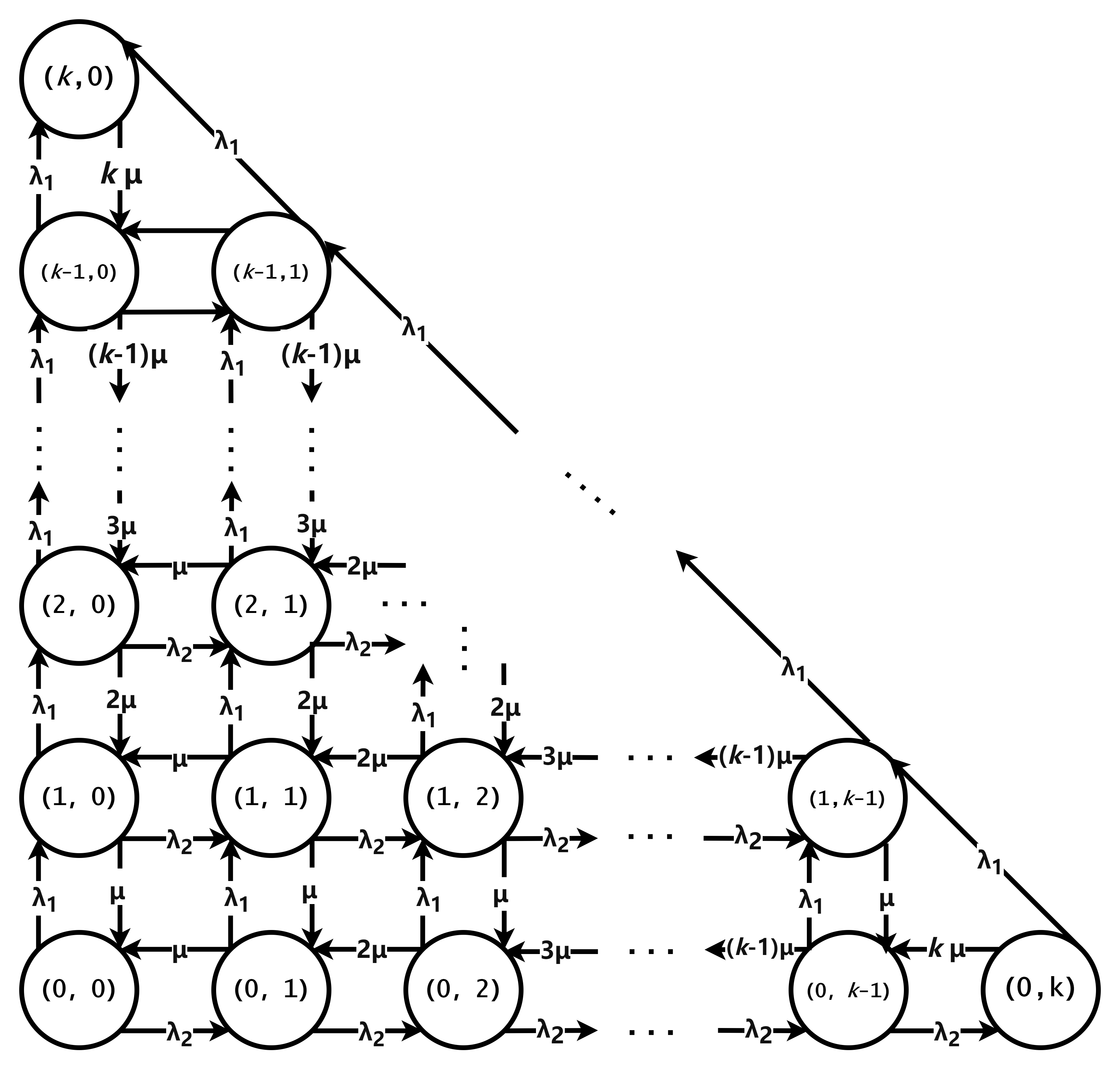}
\caption{The state transition diagram for a $k$-server loss system with 2 preemptive priorities. }\label{fig:k servers}
\end{figure}

In the case of a general number of servers ($k$), the balance equations can be constructed by considering two types of states (refer to Figure~\ref{fig:k servers}): (1) non-congestion states where at least one server is idle, and (2) congestion states where all servers are busy.

In particular, for the non-congestion states the $i,j$ satisfy\\ $0 \leq  i < k, 0 \leq  j < k , 0 \leq  i+j < k$. 

These non-congestion states obey the following global balance equations, 
\begin{equation} \label{non-congestion-general}
    \begin{split}
\displaystyle\pi_{i,j}\times\left(\lambda_1+\lambda_2+(i+j)\mu\right)=~ &\displaystyle\pi_{i-1,j}\times\lambda_1+\displaystyle\pi_{i,j-1}\times\lambda_2\\&+\displaystyle\pi_{i+1,j}\times(i+1)\mu\\&+\pi_{i,j+1}\times(j+1)\mu.
    \end{split}
\end{equation}

The congestion states can further be classified into two types: Type 1 and Type 2.\par
\par
Type 1 congestion states refer to the states in which a new arrival will be blocked only if it is of Priority class 2. It includes all the congestion states except where $i=k$. In these states, $\{i,j\}$ satisfy
$0 \leq  i < k, 0 \leq  j \leq k ,   i+j = k$. 
These Type 1 states obey the following global balance equations. 
 
\par

\begin{equation} \label{eq12_1}
    \begin{split}
\displaystyle\pi_{i,j}\times\displaystyle(\lambda_1+(i+j)\mu)=~&\displaystyle\pi_{i-1,j}\times\displaystyle\lambda_1+\displaystyle\pi_{i,j-1}\times\displaystyle\lambda_2\\&+\displaystyle\pi_{i-1,j+1}\times\displaystyle\lambda_2.
    \end{split}
\end{equation} 
\par
There is only one Type 2 state, namely when $i=k, j=0$, in which a new arrival will always be blocked irrespective of its priority class.
The state itself obeys the following global balance equation,
\par

\begin{equation} \label{eq12_2}
\pi_{k,0}\times k\mu=\pi_{k-1,1}\times\lambda_1+\pi_{k-1,0}\times\lambda_1.
\end{equation} 

\par
Note also that $\pi_{i,j}=0$ for all cases where $i<0$ or $j < 0$, and all state probabilities should sum to $1$, that is,

\begin{equation} \label{normalization}
    \begin{split}
\displaystyle\sum_{i=0}^{k}\displaystyle\sum_{j=0}^{k-i}\pi_{i,j}=1.
    \end{split}
\end{equation} 

By assuming statistical equilibrium, we can obtain the state probabilities $\pi_{i,j}$ for all possible pairs of $i$ and $j$ based on~\eqref{non-congestion-general} to~\eqref{normalization}. 
The priority 1 traffic will only be blocked in the Type 2 state ($i=k, j=0$) where all the $k$ servers are busy serving priority 1 customers. That is, 
\begin{equation} \label{eq12_3}
P_{\text{b}}(1)=\pi_{k,0} 
=\dfrac{(A_1)^k/k!}{\displaystyle\sum_{i=0}^{k}\dfrac{A_{1}^{i}}{i!}}.
\end{equation} 

For priority 2 customers, we follow similar rationales as in the previous two subsections. That is, a priority 2 customer will be blocked on arrival if the system is in either a Type 1 or Type 2 congestion state. In addition, a priority 2 customer being served will be preempted by a newly arriving priority 1 customer if the system is in a Type 1 congestion state. We again use $P_{\text{boa}}(2)$ and $P_{\text{pre}}(2)$ to represent the probabilities of priority 2 customers being blocked on arrival and preempted after being admitted to service, respectively, and derive the overall blocking probability for them as,
\begin{equation} \label{eq13}
    \begin{split}    
     P_{\text{b}}(2) =~ & P_{\text{boa}}(2) + (1-P_{\text{boa}}(2))P_{\text{pre}}(2) \\
     =~ & \displaystyle\sum_{i=0}^{k}\pi_{i,k-i}+\frac{\lambda_1}{\lambda_2}\times\left( \displaystyle\sum_{i=0}^{k}\pi_{i,k-i}-\pi_{k,0}\right)\\
     =~ & \pi_{k,0}+(1+\frac{\lambda_1}{\lambda_2})\times \left( \displaystyle\sum_{i=0}^{k}\pi_{i,k-i}-\pi_{k,0}\right)\\
     =~ & \dfrac{A_1^k/k!}{ \displaystyle\sum_{i=0}^{k}\dfrac{A_1^i}{i!}} 
     + \left(\frac{\lambda_1+\lambda_2}{\lambda_2}\right)\\&\times\left(\dfrac{(A_1+A_2)^k/k!}{ \displaystyle\sum_{i=0}^{k} \dfrac{(A_{1}+A_2)^{i}}{i!}}-\dfrac{A_1^k/k!}{ \displaystyle\sum_{i=0}^{k}\dfrac{A_1^i}{i!}}\right)\\
     =~& \dfrac{(A_1+A_2)E_k(A_1+A_2)-A_1E_k(A_1)}{A_2}.
  \end{split}
\end{equation} 

Again, this is consistent with the results in Section~\ref{V2}.

\subsection{Extension to $p$ preemptive priorities}\label{sec:subhead1}
\begin{figure}[htpb]
\centering
\includegraphics[width=0.75\columnwidth,height=\textwidth,keepaspectratio]{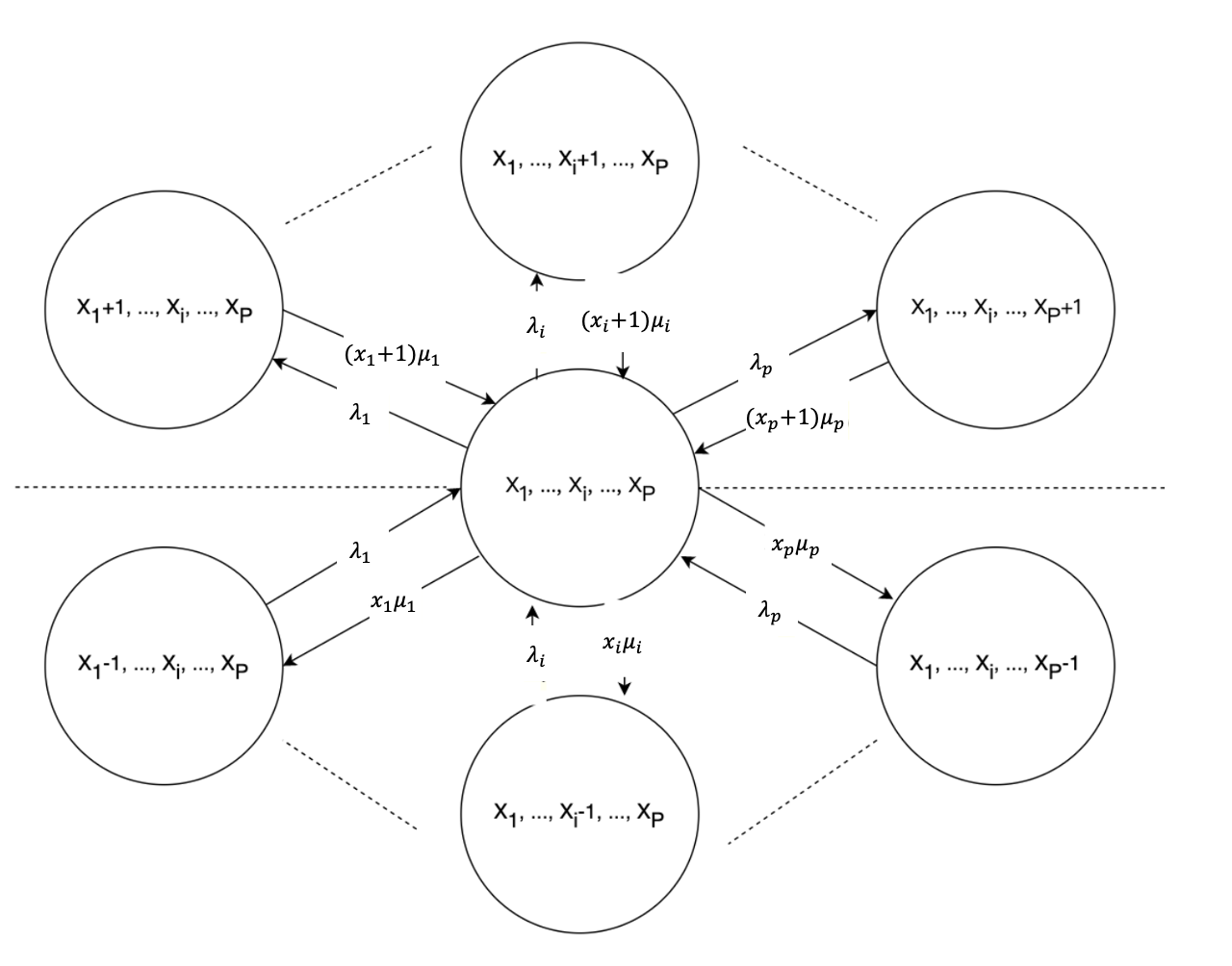}
\caption{The transition diagram of non-congestion states for a $k$-server loss system with $p$ preemptive priorities. }\label{fig:p priorities}
\end{figure}

For the M/M/$k$/$k$ model, the $k$-server loss system has a Poisson arrival process with rate $\lambda$ and exponentially distributed service time with intensity $\mu$. With $p$ priorities (the smallest number represents the highest priority), the total offered traffic is
\begin{equation} \label{eq3}
  \begin{split}
    A=\displaystyle\sum_{i=1}^{p} A_i.
  \end{split}
\end{equation}

The vector $\mathbf{X}= (...,X_i,...)$ denotes the state of the system where $X_i$ represents the number of busy servers for priority class $i$ traffic. Accordingly, $i$ and $X_i$ are bounded as follows: 
\begin{equation} \label{eq4}
  \begin{split}
    1 &\leq  i \leq p \\
    0 &\leq  X_i \leq k  \\
    0 &\leq   \displaystyle\sum_{i=1}^p X_i \leq k. \\
  \end{split}
\end{equation}

 Let $\pi(X_1..., X_i, ..., X_p)$ be the steady-state probability that there are $X_i$ busy servers of priority class $i$, for $i = 1, 2, ..., p$. Then we have:

\begin{equation} \label{eq5_2}
  \begin{split}
   \displaystyle\sum_{\begin{subarray}~X_i\in[0,k],i\in[1,p]\\ 0 \leq \sum_{i=1}^p X_i \leq k\end{subarray}}\pi(X_1..., X_i, ..., X_p)=1.
  \end{split}
\end{equation}

All the states can be divided into two types: (1) {\it non-congestion states} where for every $i = 1, 2, \ldots, p$,  
$0 \leq  X_i < k$, and 
$0 \leq   \displaystyle\sum_{i=1}^p X_i < k $ (called ``normal'' states in  \citep{Yang14}); and (2) {\it congestion states} where for every $i = 1, 2, \ldots, p$,   
$0 \leq  X_i \leq k$, and
$\displaystyle\sum_{i=1}^p X_i= k$ 
(called ``boundary'' states in \citep{Yang14,teletraffic}).


In Figure \ref{fig:p priorities}, we provide the transition diagram of the $k$-server loss system with $p$ preemptive priorities in non-congestion states (a similar figure is provided in \citep{Yang14}). This leads to the following balance equations for the non-congestion states:
\begin{equation} \label{non-congestive}
  \begin{split}
\pi(..., X_i,...) \left(\displaystyle\sum_{i=1}^{p}\lambda_i+\displaystyle\sum_{i=1}^{p}X_i\mu\right)=&
\displaystyle\sum_{i=1}^{p}\pi(... ,X_i+1,...)(X_i+1)\mu\\&+\displaystyle\sum_{i=1}^{p}\pi(..., X_i-1, ...)\lambda_i,
  \end{split}
\end{equation}
where the expression $\sum_{i=1}^{p}\pi(..., X_i+1,...)$ denotes all the possible states that can transit in one hop to  $(..., X_i,...) $ by a departure. Similarly, the expression $\sum_{i=1}^{p}(..., X_i-1, ...)$ represents all the possible states that can transit in one hop to  $(..., X_i,...) $ by an arrival. Note that $\pi(..., X_i,...) 
 = 0$ if any of the $X_i$ values is negative.
\begin{figure}[htpb]
\centering
\includegraphics[width=0.68\columnwidth,height=\textwidth,keepaspectratio]{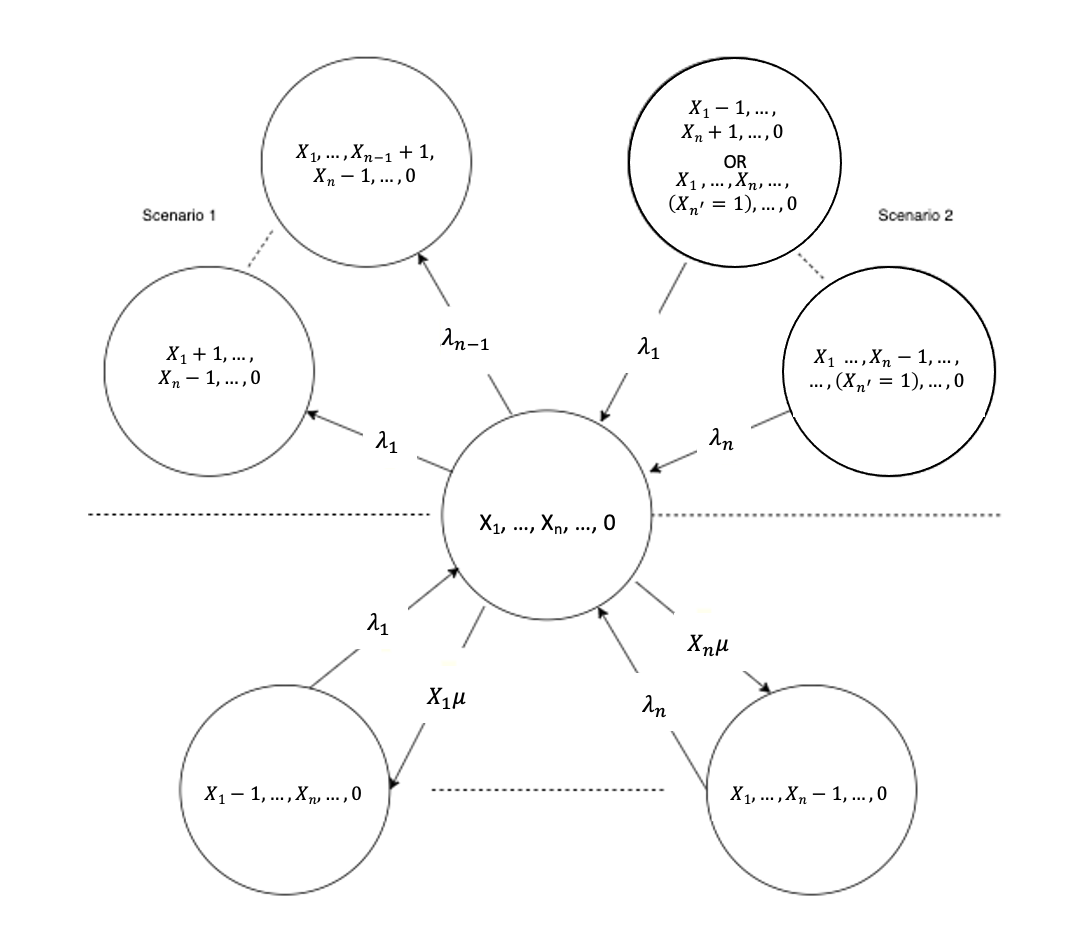}
\caption{The transition diagram of a boundary state $\mathbf{X}$, where $1 <l(\mathbf{X}) = n < p$, for a $k$-server loss system with $p$ preemptive priorities.}\label{fig:p priorities 2}
\end{figure}

Figure \ref{fig:p priorities 2} illustrates the transition diagram in boundary states. For illustration, we use $l(\mathbf{X})$ to denote the lowest priority of all customers in service in a boundary state $\mathbf{X}$. For example, $l(X_1,...,X_n,0,...,0) = n$ if $X_n > 0$ and $X_i = 0$ for all $i > n$. On the top half of the diagram, we point out two preemption scenarios for a boundary state $\mathbf{X}$ when $1 < l(\mathbf{X}) = n < p$. In scenario 1, $\mathbf{X}$ will transit to another state when a class $n$ user is preempted by an arrival of a higher class. Scenario 2, on the contrary, includes all situations where another state $\mathbf{X}' = (X'_1,...X'_p)$ transits to $\mathbf{X}$ by a single preemption. This can be further divided to two sub-scenarios: 1) if $l(\mathbf{X}') = n$, then the incoming arrival must have a higher priority than $n$, and it will preempt an $n$ priority customer in service; 2) if $l(\mathbf{X}')  = n'> n$, it implies that $X'_{n'} = 1$ and $X'_i = 0$ for all $i > n, i \neq n'$, as the state will transit to $\mathbf{X}$, where $l(\mathbf{X}) =n$, after the preemption. The incoming arrival must have a priority of $n$ or higher, and it will preempt the only $n'$ priority customer in service. 
As a consequence, equations for the boundary states when $1 < n < p$ can be written as:
\begin{equation} \label{boundary}
  \begin{split}
\pi(..., X_n,...)\left(\displaystyle\sum_{i=1}^{n-1}\lambda_i+\displaystyle\sum_{i=1}^{p}X_i\mu\right)\\=\displaystyle\sum_{i=1}^{n}\pi(..., X_i-1, ...)\lambda_i+\displaystyle\sum_{i=1}^{n-1}\pi(...,&X_i-1,...,X_{n}+1,...)\lambda_i\\+\displaystyle\sum_{i=1}^{n}\sum_{n'=n+1}^{p}\pi(...,X_i-1,&...,(X_{n'}=1),...)\lambda_i.
  \end{split}
\end{equation}

When $n=1$, it implies that $X_1 = k$ (as it is a boundary state), and there are no higher priorities that may cause a preemption. Hence:
\begin{equation} \label{n_equals_1}
\begin{array}{lcl}
\pi(k,0,...) (k\mu)&=&\pi(k-1,0,\ldots)\lambda_1\\
&&+\displaystyle\sum_{n'=2}^{p}\pi(X_1-1,...,(X_{n'}=1),...)\lambda_1.
\end{array}
\end{equation}

Respectively, when $n=p$, only customers of priority $p$ in service may be preempted by a new arrival. Hence:
\begin{equation} \label{n_equals_p}
\begin{array}{lcl}
\pi(..., X_p) \left(\displaystyle\sum_{i=1}^{p}X_i\mu+\displaystyle\sum_{i=1}^{p-1}\lambda_i\right)&\\
=\displaystyle\sum_{i=1}^{p-1}\pi(...,X_i-1,...,X_p+1)\lambda_i&+\displaystyle\sum_{i=1}^{p}\pi(..., X_i-1, ...)\lambda_i.
\end{array}
\end{equation}

Based on~\eqref{non-congestive} to~\eqref{n_equals_p}, we can obtain the probabilities of each state and, in turn, calculate the blocking probability of each priority class, given that the number of priorities $p$ is a constant.

For class 1 customers with the highest priority, blocking happens when all servers are occupied by only class 1. Therefore, the blocking probability of class 1 ($p_{1}$) can be written as $p_{1}=\pi(k,0,0,...)$.

For customers of any other priority class $n$ ($1 < n \leq p$), their blocking probabilities consist of two components: (1) the probability (denoted $P_{\text{boa}}(n)$) of being blocked because all the servers are busy and there is no user in service with a lower priority than $n$; and (2) the probability that preemption occurs  for a class $n$ customer that arrived earlier and was admitted to the system but is being preempted by an arrival of a higher priority customer (denoted $P_{\text{pre}}(n)$) which corresponds to scenario 1 shown in Figure \ref{fig:homo}. 

Let $\pi_{\text{p}n}$ be the proportion of time during which the following conditions hold. 
\begin{itemize}
\item All the $k$ servers are busy.
\item There exists at least one customer of class $n$.
\item There are no customers of priority lower than $n$ in the system.
\end{itemize}\par
During such a period of time, any arrival of a customer of priority higher than $n$ will cause preemption of a priority $n$ customer. 



Hence,

\begin{small}
\begin{equation} \label{eq6_3}
  \begin{split}
P_{\text{boa}}(n)=\dfrac{\pi \bigg( \displaystyle\sum_{i=1}^{n}X_i=k\bigg)\lambda_n}{ \lambda_n},
  \end{split}
\end{equation}

\begin{equation} \label{eq6_4}
  \begin{split}
P_{\text{pre}}(n)=\dfrac{ \displaystyle\sum_{i=1}^{n-1}\bigg(\pi_{\text{p}n}\lambda_i\bigg) }
{{(1-P_{\text{boa}}(n))}\lambda_n},
  \end{split}
\end{equation}

Define,
\begin{equation} \label{eq6_5}
  \begin{split}
P_{\text{b}}(n)=&P_{\text{boa}}(n)+{{(1-P_{\text{boa}}(n))}}P_{\text{pre}}(n)\\
=~&\dfrac{\pi\bigg( \displaystyle\sum_{i=1}^{n}X_i=k\bigg) \lambda_n }{ \lambda_n}+\dfrac{ \displaystyle\sum_{i=1}^{n-1}\bigg(\pi_{\text{p}n}\lambda_i\bigg)}{\lambda_n}
\\=~&\pi\bigg( \displaystyle\sum_{i=1}^{n}X_i=k\bigg) + \dfrac{ \displaystyle\sum_{i=1}^{n-1}\bigg(\pi_{\text{p}n}\lambda_i\bigg)}{\lambda_n}.
  \end{split}
\end{equation}
\end{small}



We generalize the model from the previous two subsections to $p$ priorities class where $A=\displaystyle\sum_{i=1}^{p} A_i$, and obtain
\begin{equation} \label{eq6_6}
\pi\bigg(\displaystyle\sum_{i=1}^{p}X_i=k\bigg)=\dfrac{A^k/k!}{\displaystyle\sum_{i=0}^{k}\dfrac{A^i}{i!}}.
\end{equation}

Similarly, by defining

\begin{equation}
\hat{A}_n = \bigg(\displaystyle\sum_{i=1}^{n}A_i\bigg),
\end{equation}

we have
\begin{equation}
P_{\text{boa}}(n)=\pi\bigg(\displaystyle\sum_{i=1}^{n}X_i=k\bigg)=\dfrac{\hat{A}^k_n/k!}{\displaystyle\sum_{j=0}^{k}\dfrac{\hat{A}^j_n}{j!}},
\end{equation}
and

\begin{equation}
\pi_{\text{p}n}=\pi\bigg(\displaystyle\sum_{i=1}^{n}X_i=k\bigg)-\pi\bigg(\displaystyle\sum_{i=1}^{n-1}X_i=k\bigg).
\end{equation}



To obtain the probability of preemption $P_{\text{pre}}(n)$, we consider the ratio of the number of preemptions that occur during an arbitrarily long period of time $L$ to the number of priority $n$ arrivals during the same period of time $L$ as follows, 

\begin{small}
\begin{equation} \label{eq54}
P_{\text{pre}}(n)=\dfrac{\Bigg(\pi\bigg(\displaystyle\sum_{i=1}^{n}X_i=k\bigg)-\pi\bigg(\displaystyle\sum_{i=1}^{n-1}X_i=k\bigg)\Bigg)\displaystyle\sum_{i=1}^{n-1}\lambda_i L}{{{(1-P_{\text{boa}}(n))}}\lambda_n L}.
\end{equation}
\end{small}

This leads to 
\begin{small}
\begin{equation} \label{eq55}
  \begin{split}
P_{\text{pre}}(n)=\dfrac{\left(\dfrac{\hat{A}_n^k/k!}{\displaystyle\sum_{m=0}^{k}\dfrac{\hat{A}^m_n}{m!}}-\dfrac{\hat{A}_{n-1}^{k}/k!}{\displaystyle\sum_{m=0}^{k}\dfrac{\hat{A}^m_{n-1}}{m!}}\right)\displaystyle\sum_{i=1}^{n-1}\lambda_i }{{{(1-P_{\text{boa}}(n))}}\lambda_n}.
  \end{split}
\end{equation}
\end{small}

Then, we obtain the following result for the overall blocking probability of priority $n$ customers. 
\begin{small}
\begin{equation} \label{eq56}
  \begin{array}{lcl}
P_{\text{b}}(n)&=&  P_{\text{boa}}(n)+{(1-P_{\text{boa}}(n))}P_{\text{pre}}(n)\\
&=&\dfrac{\hat{A}^k_n/k!}{\displaystyle\sum_{m=0}^{k}\dfrac{\hat{A}^m_n}{m!}}\\
 &&+ \dfrac{\displaystyle\sum_{i=1}^{n-1}\lambda_i \left(\dfrac{\hat{A}_n^k/k!}{\displaystyle\sum_{m=0}^{k}\dfrac{\hat{A}^m_n}{m!}}-\dfrac{\hat{A}_{n-1}^{k}/k!}{\displaystyle\sum_{m=0}^{k}\dfrac{\hat{A}^m_{n-1}}{m!}}\right)}{\lambda_n}.
\end{array}
\end{equation}
\end{small}

\subsection{Consistency with the results in Section~\ref{V2}}
As expected, the blocking probability results obtained in this paper directly from the multi-dimensional steady-state equations are consistent with the results described in Section~\ref{V2} based on \citep{V02}.

By rearranging the terms in~\eqref{eq56}, we have

\begin{small}
\begin{equation} \label{pb_general}
\begin{split}
P_{\text{b}}(n)
&=\left(\dfrac{\displaystyle\sum_{i=1}^{n}\lambda_i}{\lambda_n}\right)\dfrac{\hat{A}^k_n/k!} {\displaystyle\sum_{m=0}^{k}\dfrac{\hat{A}^m_n}{m!}}-\dfrac{\displaystyle\sum_{i=1}^{n-1}\lambda_i}{\lambda_n} \left(\dfrac{\hat{A}_{n-1}^{k}/k!}{\displaystyle\sum_{m=0}^{k}\dfrac{\hat{A}^m_{n-1}}{m!}}\right).
  \end{split}
\end{equation}
\end{small}

Recall that the service time of all customers (of all priority classes) is exponentially distributed with parameter $\mu$. By multiplying the numerators and denominators of both terms in~\eqref{pb_general} by $1/\mu$, we obtain

\begin{small}
\begin{equation} \label{eq58}
P_{\text{b}}(n)
=\left(\dfrac{\displaystyle\sum_{i=1}^{n}\dfrac{\lambda_i}{\mu}}{\dfrac{\lambda_n}{\mu}}\right)\left(\dfrac{\dfrac{\hat{A}^k_n}{k!}} {\displaystyle\sum_{m=0}^{k}\dfrac{\hat{A}^m_n}{m!}}\right)-\left(\dfrac{\displaystyle\sum_{i=1}^{n-1}\dfrac{\lambda_i}{\mu} }{\dfrac{\lambda_n}{\mu}}\right)\left(\dfrac{\dfrac{\hat{A}_{n-1}^{k}}{k!}}{\displaystyle\sum_{m=0}^{k}\dfrac{\hat{A}^m_{n-1}}{m!}}\right).
\end{equation}
\end{small}

This leads to
\begin{small}
\begin{equation}
\label{bp_proved}
P_{\text{b}}(n)=\frac{\left(\hat{A}_n\right)E_k \left(\hat{A}_n \right)-\left(\hat{A}_{n-1}\right)E_k \left(\hat{A}_{n-1}\right)}{A_n},
\end{equation}
\end{small}
which is consistent with~\eqref{lossdiff} and~\eqref{bp} as in Section~\ref{V2} and~\citep{V02}.

\section{Sensitivity to holding time distribution}\label{sec:sens}
As discussed in Section~\ref{V2},  the stationary distribution of the number of priority 1 customers in the system is insensitive to the shape of the holding time distribution. However, this insensitivity property does not apply to the preempted customers, that is, to customers of priority $i$ for  $i>1$. This can be explained as follows. Consider, for example, the case where all the service times are deterministic (fixed). When a higher-priority customer preempts a low-priority customer, the remaining service time of the preempted call is shorter than that of the higher-priority customer. In this case, the preemption will increase the load in the system, and since the displacement/replacement argument relies on the exponential service times assumption where the remaining service times of the preempted call and the new high-priority call are the same, we can expect that for the case of deterministic service times, assuming the same mean service time, we will expect to have a higher blocking probability for the preempted priorities traffic than for the case of exponential service times. On the other hand, if the service times have a higher variance than that of the exponential, the opposite will occur. Namely, the remaining service time of the preempted call is longer than that of the higher-priority customer. In this case, the preemption will decrease the load in the system, so we will expect to have a lower blocking probability for the preempted priorities traffic than for the case of exponential service times. 

Interestingly, high service time  variance improved system performance for low-priority customers. This contradicts the ``normal'' effects of queueing systems where larger variance adversely affects performance, but we have to remember that when the variance is large, mostly the longer (low priority) jobs are the ones that are being preempted by high priority traffic and this leads to the overall reduction of blocking probability.

\begin{figure*}[t]
\centering
\subfigure[]{\includegraphics[width=0.49\linewidth]{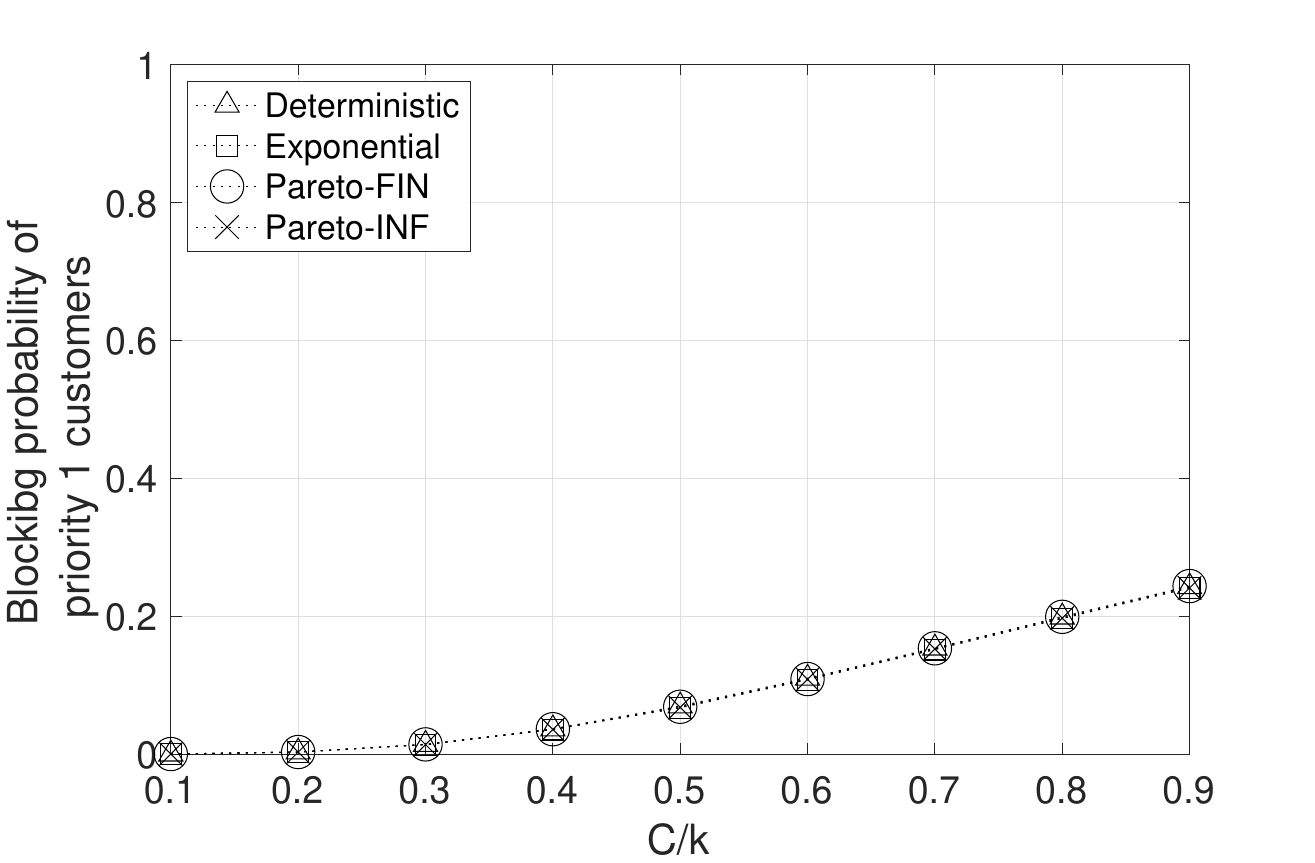}\label{fig:homo-1}}
\subfigure[]{\includegraphics[width=0.49\linewidth]{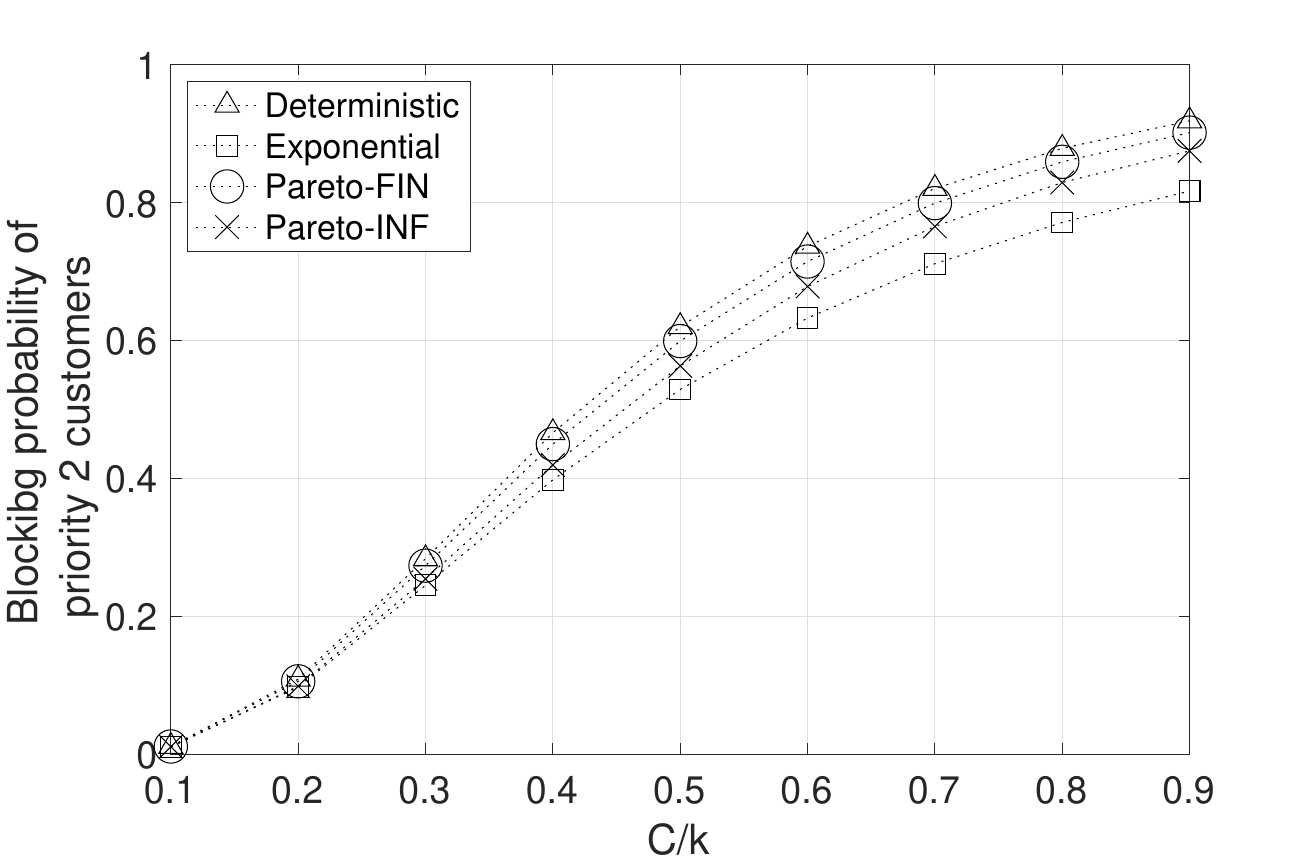}\label{fig:homo-2}}
\subfigure[]{\includegraphics[width=0.49\linewidth]{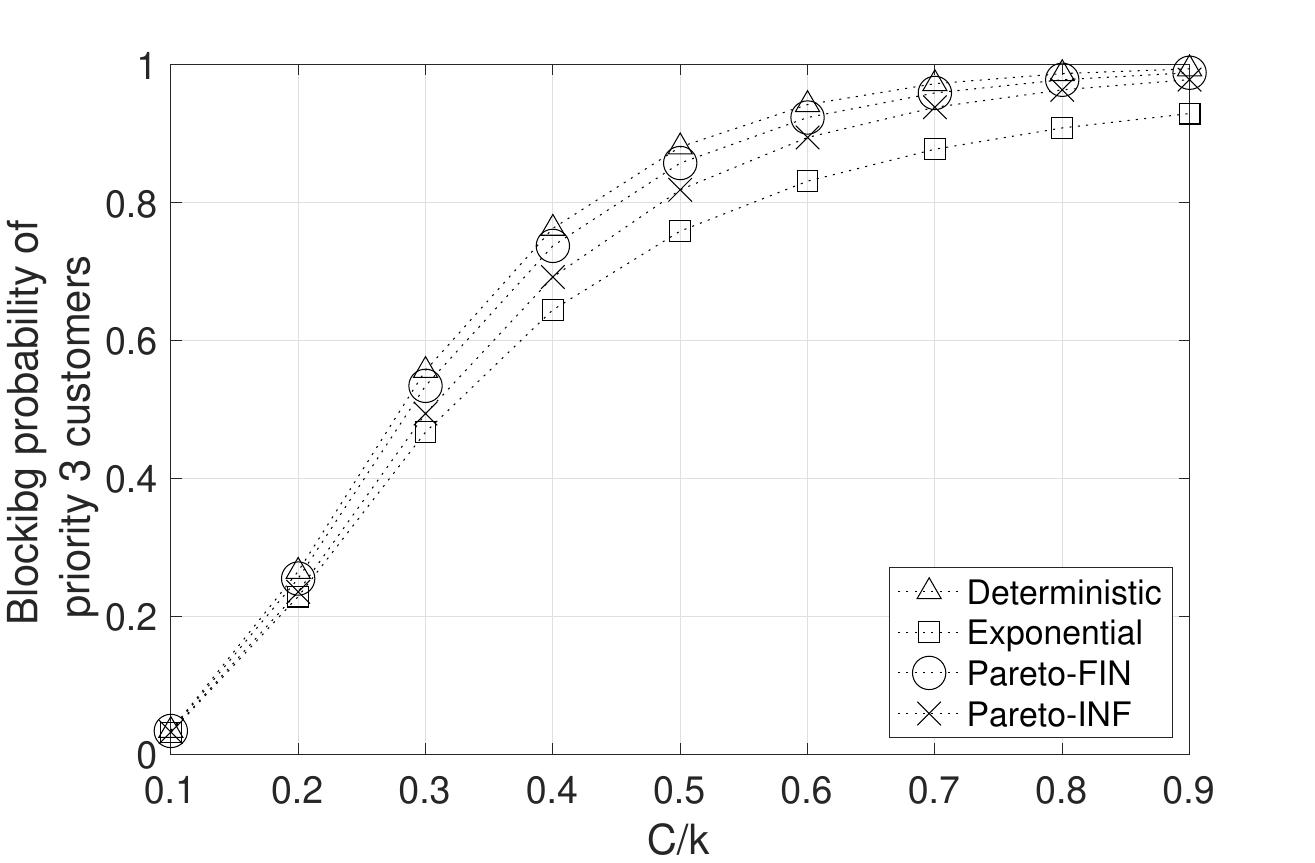}\label{fig:homo-3}}
\subfigure[]{\includegraphics[width=0.49\linewidth]{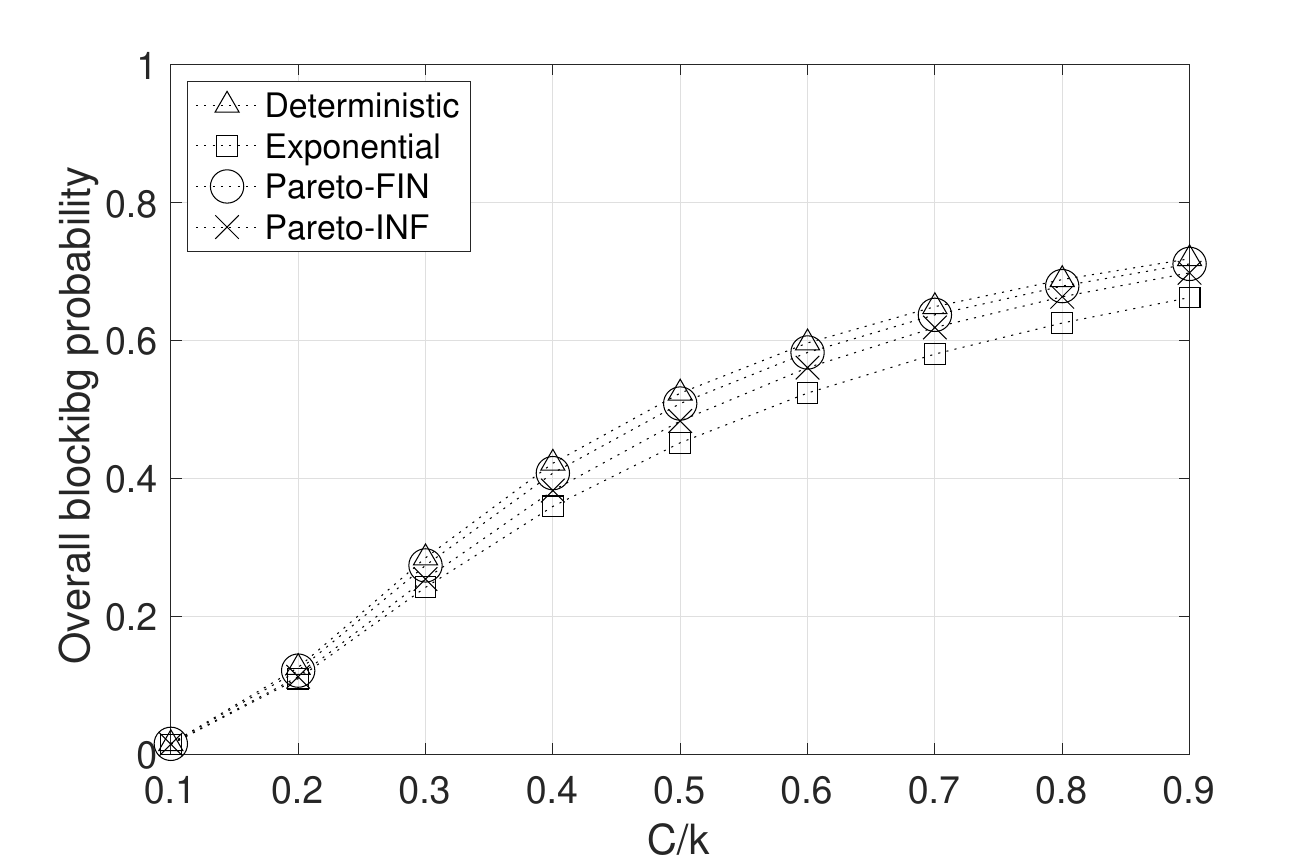}\label{fig:homo-total}}
\caption{Blocking probabilities of customers with different priorities and the blocking probability of all customers.} \label{fig:homo}
\end{figure*}

Note that this is consistent with similar comments made in \cite{klimenoklack} based on the numerical results for the loss probability of the $BMAP/PH/N/N$ queue.

We will numerically demonstrate these effects by simulations where we consider four cases involving four holding time distributions: deterministic, exponential, and two cases of Pareto distributions. In particular, we will demonstrate that the higher the variance, the lower the blocking probability of the lower-priority customers. 

We consider an M/M/$k$/$k$ system with $p=3$ priorities of customers and $k=5$ servers where the offered traffic for each priority is fixed and is denoted by $C$. That is, $A_i = C$ for $i=1,2,3$. Thus, $A=A_1+A_2+A_3=3C$.

In Figure~\ref{fig:homo}, we present our simulation results for the blocking probabilities of customers with different priorities ($p=1,2,3$) and the overall blocking probability of all customers, as a function of the ratio $C/k=C/5$. In our simulation results, the  $95\%$ confidence intervals based on the Student $t$-distribution are always maintained within $\pm 3\%$ of the observed mean.

In all the cases we consider that involve all four holding time distributions: deterministic, exponential, and the two Pareto distributions have unit mean. That is, the mean holding time is equal to one, so for the deterministic case, the service times are always equal exactly to one and for the exponential case, the service times are exponentially distributed with parameter $\mu=1$. 
We provide the blocking probabilities in a range of $C$ values between 0 and $k=5$, so that the ratio $C/k$ is between 0 and 1.
As we set $\mu=1$, for any given $C$ value, the arrival rate is given by $\lambda=C\mu=C$ for all four cases. 

For the two Pareto cases, we again consider unit mean. In the first of these two cases, the shape parameter of the Pareto distribution of the service time is set at  $\gamma = 2.001$, and for the second case, we set $\gamma = 1.98$. It is known that for a Pareto random variable, the variance is infinite for $0< \gamma \leq 2$. Accordingly, we denote the first case, for which the variance is finite (with $\gamma = 2.001$),  as Pareto-FIN and the second case (with $\gamma = 1.98$, where the variance is infinite) as Pareto Pareto-INF. Both cases represent distributions with very large variances. 

Also, note that as the mean of the service time is set to be equal to one, the scale parameter denoted $\delta$ for each of the two cases can be obtained by $$\delta = \frac{\gamma-1}{\gamma}. $$ 
Accordingly, the values of the Pareto scale parameters are set at $\delta = 0.50025$ and $\delta = 0.495$ for Pareto-FIN and Pareto-INF, respectively. 

In Figure~\ref{fig:homo-1}, we demonstrate that the blocking probabilities of priority 1 customers are insensitive to the shape of the tested holding time distributions; while, in Figures~\ref{fig:homo-2} and \ref{fig:homo-3}, it is demonstrated that the blocking probabilities of other customers are clearly different for different distributions and that high  variance of the service time distribution leads to lower blocking probabilities for the lower priority customers. 
In consistency with the sensitivity of the preempted customers, the curves of the overall blocking probability presented in Figure~\ref{fig:homo-total} are also different from each other, and again we observe that a larger variance of the service time reduces the blocking probability.  
These numerical results are sufficient to demonstrate the sensitivity of preempted customers to the shape of their holding time distribution.

\section{Conclusion}\label{sec:conc}
We have established the consistency between the global balance (steady state) equations for a loss system with preemptive priorities and a result for the various blocking probabilities of the different classes of customers based on traffic loss arguments. 
This has been done by deriving this known result directly from the global balance equations of the relevant multidimensional Markov chain. 





This has been achieved by observing that the blocking probability of a customer of any other priority class consists of two components: (1) the probability of being blocked because all the servers are busy and there is no user in service with a lower priority, and (2) the probability that the customer has been admitted, but after admission, it has been preempted  by an arrival of a higher priority customer. Notice that for the top priority class, the second component is equal to zero. After these two probabilities (components) have been derived, we have obtained the steady-state blocking probabilities of all customers, and we have observed the consistency with the previously obtained results based on traffic loss arguments.

We have also provided explanations and demonstrated by simulations that except for the blocking probability of the highest priority customers, the blocking probabilities of the other customers are sensitive to the holding time distributions and that a higher variance of the service time distributions reduces the blocking probabilities of the lower priority customers. 



\section*{Acknowledgment}

The authors would like to thank Prof. Peter Taylor for his helpful comments on the paper.
The authors are also thankful to the anonymous reviewers for their helpful comments, corrections, and suggestions.


\bibliographystyle{plos2015}
\bibliography{references_clean}

\end{document}